\documentclass[12pt,oneside,english]{amsart}
\usepackage[latin1]{inputenc}
\usepackage{amssymb}

\theoremstyle{plain}
\newtheorem{theorem}{Theorem}

\newtheorem{lemma}[theorem]{Lemma}

\theoremstyle{definition}

\theoremstyle{remark}
\newtheorem{remark}[theorem]{Remark}

\usepackage{babel}

\makeatother
\begin{document}
\baselineskip=17pt

\title[Digit functions with large periods]
{A family of digit functions with large periods}

\author{Vladimir Shevelev}
\address{Department of Mathematics \\Ben-Gurion University of the
 Negev\\Beer-Sheva 84105, Israel. e-mail: shevelev@bgu.ac.il}
 \author{Peter J. C. Moses}
 \address{United Kingdom. e-mail: mows@mopar.freeserve.co.uk}


\begin{abstract}
For odd $n\geq3,$ we consider a general hypothetical identity for the differences
$S_{n,\;0}(x)$ of multiples of $n$ with even and odd digit sums in the base $n-1$
in interval $[0,x),$ which we prove in the cases $n=3$ and $n=5$ and empirically
 confirm for some other $n.$ We give a verification algorithm for this identity for
 any odd $n.$ The hypothetical identity allows to give a general recursion for
 $S_{n,\;0}(x)$ for every integer $x$ depending on the residue of $x$ modulo
  $p(n)=2n(n-1)^{n-1},$ such that $p(3)=24, \;p(5)=2560,\;p(7)=653184,\; etc.$
\end{abstract}

\maketitle
\section{Introduction}
For $x\in \mathbb{N}$ and $n\geq 3$, denote by $S_{n}(x)$ the sum

\begin{equation}\label {1}
S_{n,\;j}(x)=\sum_{0\leq r < x: \;\;r\equiv j\pmod n}(-1)^{s_{n-1}(r)},
\end{equation}
where $s_{n-1}(r)$ is the digit sum of $r$ in base $n-1.$\newline
Note that, in particular, $S_{3,\;0}(x)$ equals the difference between the numbers
 of multiples of 3 with even and odd binary digit sums (or multiples of
 3 from sequences A001969 and A000069 in \cite{7}) in interval $[0,x).$

Leo Moser (cf. \cite{3}, Introduction) conjectured that always

\begin{equation}\label{2}
S_{3,\;0}(x)>0.
\end{equation}

 Newman \cite{3} proved this conjecture. Moreover, he obtained
the inequalities

\begin{equation}\label{3}
\frac{1}{20} < S_{3,\;0}(x)x^{-\lambda}< 5,
\end{equation}
where
\begin{equation}\label{4}
\lambda=\frac{\ln 3}{\ln 4}=0.792481...\;.
\end{equation}
In connection with this, the qualitative result
 (\ref{2}) we call a weak Newman phenomenon (or Moser-Newman phenomenon),
 while an estimating result of the form (\ref{3}) we call a strong Newman
 phenomenon.\newline
\indent In 1983, Coquet \cite{1} studied a very complicated continuous and
nowhere differentiable fractal
function $F(x)$ with period 1 for
which

\begin{equation}\label{5}
S_{3,\;0}(3x)=x^\lambda F\left(\frac{\ln x}{\ln
4}\right)+\frac{\eta(x)}{3},
\end{equation}

where

\begin{equation}\label{6}
\eta(x)=\begin{cases} 0,\;\; if \; x \;\; is\;\; even,\\
(-1)^{s_2(3x-1)}, \;\; if \;\; x \;\; is\;\; odd.\end{cases}
\end{equation}

He obtained that

\begin{equation}\label{7}
\limsup_{x\rightarrow
\infty, \;x\in\mathbb{N}}S_{3,\;0}(3x)x^{-\lambda}=\frac{55}{3}\left(\frac{3}{65}\right)
^\lambda=1.601958421\ldots\;,
\end{equation}

\begin{equation}\label{8}
\liminf_{x\rightarrow\infty, \;x\in\mathbb{N}}S_{3,\;0}(3x)x^{-\lambda}=\frac{2\sqrt{3}}{3}
=1.154700538\ldots\;.
\end{equation}

In 2007, Shevelev \cite{4} gave an elementary proof of Coquet's formulas
(\ref{7})-(\ref{8}) and his sharp estimates in the form
\begin{equation}\label{9}
\frac{2\sqrt{3}}{3}x^\lambda\leq
S_{3\;0}(3x,\;0)\leq\frac{55}{3}\left(\frac {3}{65}\right)^\lambda
x^\lambda,\;\;x\in\mathbb{N}.
\end{equation}
In \cite{4} it was found the following simple identity
\begin{equation}\label{10}
S_{3,\;0}(4x)= 3S_{3,\;0}(x),\;where\; x\; is \; even.
\end{equation}
Since in the left hand side of (\ref{10}) the argument $4x\equiv 0 \pmod 8$
then (\ref{10}) is not a recursion for evaluation of $S_{3,\;0}(x).$
However, in the same work Shevelev found the following recursion
for fast calculation of $S_{3,\;0}(x):$
\begin{equation}\label{11}
S_{3,0}(x)=3S_{3,0}\left(\left\lfloor \frac x
4\right\rfloor\right)+\nu(x),
\end{equation}

where

\begin{equation}\label{12}
\nu(x)=\begin{cases} 0,\;if \;x\equiv 0,7,8,9,16,17,18,22,23\;\pmod
{24};\\ (-1)^{s_2(x)},\; if \; x\equiv 3,4,10,12,20 \; \pmod {24};\\
(-1)^{s_2(x)+1},\; if \; x\equiv 1,2,5,6,11,19,21\; \pmod {24};\\
2(-1)^{s_2(x)},\; if \; x\equiv 15\;\pmod {24};\\
2(-1)^{s_2(x)+1},\; if \; x\equiv 13,14\;\pmod {24}.\end{cases}
\end{equation}

 In 2008, Drmota and Stoll \cite{2} proved a generalized weak Newman phenomenon,
showing that (\ref{2}) is valid for $S_{n,\;0}(x)$ for every $n\geq3,$
 at least beginning with $x\geq x_0(n).$ A year before, Shevelev \cite{5} proved
 a strong form
 \newpage
 of this generalization, but yet only in "full"
   intervals of the form $[0, (n-1)^{2p}).$ Recently Shevelev and Moses \cite{6} in the case
of odd $n\geq3$ and $p\geq\frac{n-1}{2}$ found the relation
 \begin{equation}\label{13}
\sum_{k=0}^{\frac{n-1}{2}}(-1)^k\binom{n}{2k}S_{n,\;0}((n-1)^{2p-2k})=
\begin{cases} 0,\;\; if \; p\geq\frac{n+1}{2},
\\(-1)^n, \;\; if \;\;p=\frac{n-1}{2} .\end{cases}
\end{equation}
In the case of $p=\frac{n-1}{2},$ (\ref{13}) could be rewrite in the form
\begin{equation}\label{14}
\sum_{j=0}^{\frac{n-1}{2}}(-1)^j\binom{n}{2j+1}S_{n,\;0}((n-1)^{2j})=1.
\end{equation}
Numerous experiments show that, most likely, the following more general
relation takes place:
$$\sum_{j=0}^{\frac{n-1}{2}}(-1)^j\binom{n}{2j+1}S_{n,\;0}((n-1)^{2j}x)=$$
\begin{equation}\label{15}
\sum_{j=0}^{n-1}S_{n,\;j}(x),\;x\geq1,\; n\equiv1\pmod2.
\end{equation}
In particular, we verified (\ref{15}) for $n=3,5,7,...,35$ and $1\leq x\leq1000.$
It is clear that (\ref{14}) is a special case of (\ref{15}) for $x=1,$ since
\begin{equation}\label{16}
 S_{n,\;j}(1)=\begin{cases} 1,\;\; if \; j=0,
\\0, \;\; if \;\;1\leq j\leq n-1.\end{cases}
\end{equation}
Below we show that (\ref{15})  allows with the uniform positions
to find a recursion for $S_{n,\;0}(x)$ for every odd $n\geq3.$ In the two first sections
we prove identity (\ref{15}) in cases $n=3$ and $n=5.$ In Section 4 we give a
general verification algorithm for the identity (\ref{15}) which
 allows to prove the identity (\ref{15}) for $n=7,9,..., etc. $
  In Section 5 we give a simplification of the conjectural equality (\ref{15}).
  In Section 6 we prove the recursion in case $n=3$ and in Section 7 we give the
  recursion in case $n=5.$ After these sections, in supposition that (\ref{15}) is
  true, it will be clear how to find the further recursions for odd $n\geq7.$

 \section{The identity in case $n=3$}
 Note that, by (\ref{1}),

 $$S_{3,\;j}(x)=\sum_{0\leq r < x: \;\;r\equiv j\pmod 3}(-1)^{s_{2}(r)}$$
 which yields that
  \begin{equation}\label{17}
\sum_{0\leq r < 2x: \;\;r\equiv 2j\pmod 6}(-1)^{s_{2}(r)},\;j=0,1,2.
 \end{equation}
 \newpage
 On the other hand,
  $$S_{3,\;j}(2x)=\sum_{0\leq r < 2x: \;\;r\equiv j\pmod 6}(-1)^{s_{2}(r)}+$$

\begin{equation}\label{18}
\sum_{0\leq r < 2x: \;\;r\equiv j+3\pmod 6}(-1)^{s_{2}(r)},\;j=0,1,2.
\end{equation}
 Using (\ref{18}), for $j=0,1,2,$ we consecutively find

 $$S_{3,\;0}(2x)=\sum_{0\leq r < 2x: \;\;r\equiv 0\pmod 6}(-1)^{s_{2}(r)}-$$
\begin{equation}\label{19}
\sum_{0\leq r < 2x: \;\;r\equiv 2\pmod 6}(-1)^{s_{2}(r)},
\end{equation}
$$S_{3,\;1}(2x)=-\sum_{0\leq r < 2x: \;\;r\equiv 0\pmod 6}(-1)^{s_{2}(r)}+$$
\begin{equation}\label{20}
\sum_{0\leq r < 2x: \;\;r\equiv 4\pmod 6}(-1)^{s_{2}(r)},
\end{equation}
$$S_{3,\;2}(2x)=\sum_{0\leq r < 2x: \;\;r\equiv 2\pmod 6}(-1)^{s_{2}(r)}-$$
\begin{equation}\label{21}
\sum_{0\leq r < 2x: \;\;r\equiv 4\pmod 6}(-1)^{s_{2}(r)}.
\end{equation}
Now the application of (\ref{17}) to (\ref{19})-(\ref{21}) yields the relations
\begin{equation}\label{22}
S_{3,\;0}(2x)=S_{3,\;0}(x)-S_{3,\;1}(x),
 \end{equation}
 \begin{equation}\label{23}
S_{3,\;1}(2x)=-S_{3,\;0}(x)+S_{3,\;2}(x),
 \end{equation}
 \begin{equation}\label{24}
S_{3,\;2}(2x)=S_{3,\;1}(x)-S_{3,\;2}(x).
 \end{equation}
 For $n=3,$ the left hand side of (\ref{15}) is $3S_{3,\;0}(x)-S_{3,\;0}(4x)$ and,
  using (\ref{22})-({24}), we have
  $$3S_{3,\;0}(x)-S_{3,\;0}(4x)=3S_{3,\;0}(x)-S_{3,\;0}(2x)+S_{3,\;1}(2x)=$$ $$3S_{3,\;0}(x)-S_{3,\;0}(x)+S_{3,\;1}(x)-S_{3,\;0}(x)+S_{3,\;2}(x)=$$
  $$S_{3,\;0}(x)+S_{3,\;1}(x)+S_{3,\;2}(x) $$
  which proves (\ref{15}) in the case $n=3.$
  \newpage
  \section{The identity in case $n=5$}
  In the same way, instead of (\ref{22})-(\ref{24}), we find the following relations
 \begin{equation}\label{25}
S_{5,\;0}(4x)=S_{5,\;0}(x)-S_{5,\;1}(x)+S_{5,\;2}(x)-S_{5,\;3}(x),
 \end{equation}
 \begin{equation}\label{26}
S_{5,\;1}(4x)=-S_{5,\;0}(x)+S_{5,\;1}(x)-S_{5,\;2}(x)+S_{5,\;4}(x),
 \end{equation}
 \begin{equation}\label{27}
S_{5,\;2}(4x)=S_{5,\;0}(x)-S_{5,\;1}(x)+S_{5,\;3}(x)-S_{5,\;4}(x),
 \end{equation}
 \begin{equation}\label{28}
S_{5,\;3}(4x)=-S_{5,\;0}(x)+S_{5,\;2}(x)-S_{5,\;3}(x)+S_{5,\;4}(x),
 \end{equation}
\begin{equation}\label{29}
S_{5,\;4}(4x)=S_{5,\;1}(x)-S_{5,\;2}(x)+S_{5,\;3}(x)-S_{5,\;4}(x).
 \end{equation}
For $n=5,$ the left hand side of (\ref{15}) is
\begin{equation}\label{30}
 5S_{5,\;0}(x)-10S_{5,\;0}(16x)+S_{5,\;0}(256x).
\end{equation}
Using (\ref{25})-(\ref{29}), we easily find
 \begin{equation}\label{31}
S_{5,\;0}(16x)=4S_{5,\;0}(x)-3S_{5,\;1}(x)+S_{5,\;2}(x)+S_{5,\;3}(x)-3S_{5,\;4}(x),
 \end{equation}
 \begin{equation}\label{32}
S_{5,\;1}(16x)=-3S_{5,\;0}(x)+4S_{5,\;1}(x)-3S_{5,\;2}(x)+S_{5,\;3}(x)+S_{5,\;4}(x),
 \end{equation}
 \begin{equation}\label{33}
S_{5,\;2}(16x)=S_{5,\;0}(x)-3S_{5,\;1}(x)+4S_{5,\;2}(x)-3S_{5,\;3}(x)+S_{5,\;4}(x),
 \end{equation}
 \begin{equation}\label{34}
S_{5,\;3}(16x)=S_{5,\;0}(x)+S_{5,\;1}(x)-3S_{5,\;2}(x)+4S_{5,\;3}(x)-3S_{5,\;4}(x),
 \end{equation}
 \begin{equation}\label{35}
S_{5,\;4}(16x)=-3S_{5,\;0}(x)+S_{5,\;1}(x)+S_{5,\;2}(x)-3S_{5,\;3}(x)+4S_{5,\;4}(x).
 \end{equation}
 Now using (\ref{31})-(\ref{35}), we find

$$S_{5,\;0}(256x)=36S_{5,\;0}(x)-29S_{5,\;1}(x)+$$
 \begin{equation}\label{36}
11S_{5,\;2}(x)+11S_{5,\;3}(x)-29S_{5,\;4}(x).
 \end{equation}
 Finally, for the expression (\ref{30}), using (\ref{31}) and (\ref{36}), we have

$$5S_{5,\;0}(x)-10S_{5,\;0}(16x)+S_{5,\;0}(256x)=$$
\begin{equation}\label{37}
S_{5,\;0}(x)+S_{5,\;1}(x)+S_{5,\;2}(x)+S_{5,\;3}(x)+S_{5,\;4}(x).
 \end{equation}
 It is the identity (\ref{15}) in the case $n=5.$
 \newpage
  \section{General problem}
  Quite analogously to systems (\ref{22})-(\ref{24}), (\ref{25})-(\ref{29})
   we can write the system for any $n\geq3.$ For odd $n,$ we have
   $$S_{n,\;0}((n-1)x)=S_{n,\;0}(x)-S_{n,\;1}(x)+...+S_{n,\;n-3}(x)-S_{n,\;n-2}(x),$$
   $$S_{n,\;1}((n-1)x)=-S_{n,\;0}(x)+S_{n,\;1}(x)-...-S_{n,\;n-3}(x)+S_{n,\;n-1}(x),$$
   $$S_{n,\;2}((n-1)x)=S_{n,\;0}(x)-S_{n,\;1}(x)+...-S_{n,\;n-4}(x)+S_{n,\;n-2}(x)
   -S_{n,\;n-1}(x),$$
   $$................................... $$
   $$S_{n,\;n-2}((n-1)x)=-S_{n,\;0}(x)+S_{n,\;2}(x)-...-S_{n,\;n-2}(x)+S_{n,\;n-1}(x),
   $$
   \begin{equation}\label{38}
   S_{n,\;n-1}((n-1)x)=S_{n,\;1}(x)-S_{n,\;2}(x)+...+S_{n,\;n-1}(x).
   \end{equation}
   It is easy to see that the right hand side of the $i$-th equality
    for $S_{n,\;i}((n-1)x), \;i=0,1,...,n-1,$ of the system (\ref{38}) satisfies
   the rules: 1) the signs alternate, beginning with $(-)^i;$\;
   2) there is no summand $S_{n,\;n-1-i}(x).$ Using, as usual, the convention
    $\sum_a^b=0,$ if $b<a,$ one can write
   the system (\ref{38}) in the form
   \begin{equation}\label{39}
   (-1)^iS_{n,\;i}((n-1)x))=\sum_{j=0}^{n-i-2}(-1)^{j}S_{n\;j}(x)-
   \sum_{j=n-i}^{n-1}(-1)^{j}S_{n\;j}(x).
   \end{equation}
   Thus the general problem is to prove that (\ref{39}) yields (\ref{15}).
  \section{A simplification of the conjecture}
  Note that in the sum $\sum_{j=0}^{n-1}S_{n,\;j}(x)$ the index of
  summing $j$ runs all residues modulo $n.$ Therefore, we have
  $$\sum_{j=0}^{n-1}S_{n,\;j}(x)=S_{1,\;0}(x)=\sum_{0\leq i<x}(-1)^{s_{n-1}(i)}=$$
  \begin{equation}\label{40}
  \begin{cases} 0,\;\; if \; x \;\; is\;\; even,\\
(-1)^{s_{n-1}(x-1)}, \;\; if \;\; x \;\; is\;\; odd.\end{cases}
 \end{equation}
 Thus the conjectural relation (\ref{15}) is equivalent to the equality
$$\sum_{j=0}^{\frac{n-1}{2}}(-1)^j\binom{n}{2j+1}S_{n,\;0}((n-1)^{2j}x)=$$
 \begin{equation}\label{41}
\begin{cases} 0,\;\; if \; x \;\; is\;\; even,\\
(-1)^{s_{n-1}(x-1)}, \;\; if \;\; x \;\; is\;\; odd.\end{cases}
\end{equation}
In particular, for $x=1,$ we again have (\ref{14}).
 Note that (\ref{41}) means that its left hand side taken with sign
 $(-1)^{s_{n-1}(x-1)}$ is periodic with period 2:
 \newpage
$$(-1)^{s_{n-1}(x-1)}\sum_{j=0}^{\frac{n-1}{2}}(-1)^j\binom{n}{2j+1}
S_{n,\;0}((n-1)^{2j}x)=$$
\begin{equation}\label{42}
\begin{cases} 0,\;\; if \; x \;\; is\;\; even,\\
1, \;\; if \;\; x \;\; is\;\; odd.\end{cases}
\end{equation}
\section{Recursion for $S_{3,0}(x)$}
Here we prove (\ref{11})-(\ref{12}). Let us write (\ref{42}) for $n=3$ and
$x:=\lfloor\frac{x}{4}\rfloor.$ We have
$$(-1)^{s_{2}(\lfloor\frac{x}{4}\rfloor-1)}(3S_{3,\;0}(\lfloor\frac{x}{4}\rfloor)-
S_{3,\;0}(4\lfloor\frac{x}{4}\rfloor)= $$
\begin{equation}\label{43}
\begin{cases} 0,\;\; if \; \lfloor\frac{x}{4}\rfloor \;\; is\;\; even,\\
1, \;\; if \;\; \lfloor\frac{x}{4}\rfloor \;\; is\;\; odd.\end{cases}
\end{equation}
Note that $\lfloor\frac{x}{4}\rfloor$ is even, if
$x=0,1,2,3,8,9,10,11,...$ and odd for other integers.
Thus we obtain
\begin{lemma} \label{l1}The sequence $\{A_3(x)\}$, where
\begin{equation}\label{44}
A_3(x)=(-1)^{s_2(\lfloor \frac x
4\rfloor-1)}(3S_{3,\;0}(\lfloor\frac{x}{4}\rfloor)-
S_{3,\;0}(4\lfloor\frac{x}{4}\rfloor),
\end{equation}
is periodic with the period $8,$ such that
\begin{equation}\label{45}
 A_3(x)=\begin{cases} 0,\;\; if \; x\equiv0,1,2,3,\pmod{8},\\
1, \;\; if \;\;x\equiv4,5,6,7\pmod{8}.\end{cases}
\end{equation}
\end{lemma}
 Consider the difference
 \begin{equation}\label{46}
\Delta_3(x)=S_{3,\;0}(x)-S_{3,\;0}(4\lfloor\frac{x}{4}\rfloor).
   \end{equation}
\begin{lemma}\label{l2}
We have
\begin{equation}\label{47}
\Delta_3(x)=\begin{cases} (-1)^{s_2(x-1)},\;\; if \; x\equiv1,7\;or\;10\pmod{12}\\
 (-1)^{s_2(x-2)}, \;\; if \; x\equiv2\;or\;11\pmod{12} \\(-1)^{s_2(x-3)}, \;\; if
 \; x\equiv3 \pmod{12}\\0, otherwise. \end{cases}
 \end{equation}
\end{lemma}
\slshape Proof.\;\;\upshape Let $x=12t+j,\;j=0,1,...,11.$ Consider 3 cases.
$$a)\; j=0,1,2\;or\;3.$$
Then
$$\Delta_3(x)=S_{3,0}(12t+j)-S_{3,0}(12t)=$$

 $$\begin{cases} 0,\;\; if \; j=0,\\
(-1)^{s_2(x-j)}, \;\; if \;\; j=1,2,3.\end{cases}$$
\newpage
$$b)\; j=4,5,6\;or\;7.$$
Then
 $$\Delta_3(x)=S_{3,0}(12t+j)-S_{3,0}(12t+4)=$$

 $$\begin{cases} 0,\;\; if \; j=4,5,6,\\
(-1)^{s_2(x-1)}, \;\; if \;\; j=7.\end{cases}$$

$$c)\; j=8,9,10\;or\;11.$$
Then
 $$\Delta_3(x)=S_{3,0}(12t+j)-S_{3,0}(12t+8)=$$

  $$\begin{cases} 0,\;\; if \; j=8,9,\\
(-1)^{s_2(x-1)}, \;\; if \;\; j=10,\\(-1)^{s_2(x-2)}, \;\; if \;\; j=11\end{cases}
\;$$
and (\ref{47}) follows.\; $\blacksquare$\newline
Now from (\ref{44})-(\ref{47}) we easily deduce the following result.
\begin{theorem}\label{t3}
 \begin{equation}\label{48}
S_{3,\;0}(x)=3S_{3,\;0}(\lfloor\frac{x}{4}\rfloor)+\Delta_3(x)-(-1)^{s_2(\lfloor
 \frac {x}4\rfloor-1)}A_3(x),
\end{equation}
where $A_3(x)$ and $\Delta_3(x)$ are defined by $(\ref{45})$ and $(\ref{47})$
 respectively.
 \end{theorem}
 Formula (\ref{48}) gives a recursion for $S_{3,\;0}(x).$ Let us show that it
 coincides with the recursion (\ref{11})-(\ref{12}), i.e.,
  \begin{equation}\label{49}
\Delta_3(x)-(-1)^{s_2(\lfloor\frac {x}4\rfloor-1)}A_3(x)=\nu(x),
\end{equation}
where $\nu(x)$ is defined by (\ref{12}). This follows from the following two lemmas.
\begin{lemma}\label{l4} The sequence
\begin{equation}\label{50}
\{(-1)^{s_2(x)+s_2(\lfloor\frac {x}4\rfloor-1)}A_3(x)\}
\end{equation}
is periodic with period 8.
\end{lemma}
\slshape Proof.\;\;\upshape In cases $x\equiv i\pmod8,\; i=0,1,2,3$ the terms of the
sequence are zeros. If $x\equiv i\pmod8,\; i=4,5,6,7,$ put $x=8t+i.$ Then $A_3(x)=1$
and we have
$$(-1)^{s_2(x)+s_2(\lfloor\frac {x}{4}\rfloor-1)}=(-1)^{s_2(8t+i)+s_2(2t)}=$$

 $$(-1)^{s_2(8t+i)+s_2(8t)}=(-1)^{s_2(i)}$$
 and the lemma follows. \;$\blacksquare$\newline
 Note that period of sequence (\ref{50}) is
 \newpage
  \begin{equation}\label{51}
  \{0,0,0,0,-1,1,1,-1\}.
  \end{equation}
\begin{lemma}\label{l5} The sequence
\begin{equation}\label{52}
\{(-1)^{s_2(x)}\Delta_3(x)\}
\end{equation}
is periodic with period 12.
\end{lemma}
\slshape Proof.\;\;\upshape According to (\ref{47}), we have
$$(-1)^{s_2(x)}\Delta_3(x)=$$

\begin{equation}\label{53}
\begin{cases} (-1)^{s_2(x)+s_2(x-1)},\;\; if \; x\equiv1,7\;or\;10\pmod{12}\\
(-1)^{s_2(x)+s_2(x-2)}, \;\; if \; x\equiv2\;or\;11\pmod{12} \\(-1)^{s_2(x)+s_2(x-3)}, \;\; if
 \; x\equiv3 \pmod{12}\\0, otherwise. \end{cases}
 \end{equation}
 Let $x=12t+i,\;0\leq i\leq11.$ Let, firstly, $i=1,7,10.$ In cases $i=1$ and $i=7,$ we, evidently, have
$(-1)^{s_2(x)+s_2(x-1)}=-1,$ while in case $i=10,$ $$(-1)^{s_2(12t+10)+s_2(12t+9)}
=(-1)^{s_2(12t+1010_2)+s_2(12t+1001_2)}=1.$$
Let now $i=2,11.$ In case $i=2,$ we, evidently, have
$(-1)^{s_2(x)+s_2(x-2)}=-1$ and also in case $i=11,$ we find
$$(-1)^{s_2(12t+11)+s_2(12t+9)}=(-1)^{s_2(12t+1011_2)+s_2(12t+1001_2)}=-1;$$
finally, if $i=3,$ then, evidently, we have $(-1)^{s_2(x)+s_2(x-3)}=1.$ In other
 cases, the terms of the sequence are zeros.\;$\blacksquare$\newline
 Thus period of sequence (\ref{52}) is
  \begin{equation}\label{54}
  \{0,-1,-1,1,0,0,0,-1,0,0,1,-1\}.
  \end{equation}
 Subtracting the tripled period (\ref{51}) from the doubled period (\ref{54}), we
 obtain the period of length 24 of the left hand side of (\ref{49}) multiplied by
 $(-1)^{s_2(x)}.$ It is
 $$\{0,-1,-1,1,1,-1,-1,0,0,0,1,-1,$$
   \begin{equation}\label{55}
  1,-2,-2,2,0,0,0,-1,1,-1,0,0\}.
  \end{equation}
 It is left to note that, according to (\ref{12}), $(-1)^{s_2(x)}\nu(x)$ is
 periodic with the same period.\; $\blacksquare$
\newpage
 \section{On recursion for $S_{n,0}(x)$}
Let (\ref{42}) be true. Let us write (\ref{42}) for
$x:=\lfloor\frac{x}{(n-1)^{n-1}}\rfloor.$ We have
$$(-1)^{s_{n-1}(\lfloor\frac{x}{(n-1)^{n-1}}\rfloor-1)}
((-1)^{\frac{n-1}{2}}S_{n,\;0}((n-1)^{n-1}\lfloor\frac{x}{(n-1)^{n-1}}\rfloor)+$$ $$\sum_{j=0}^{\frac{n-3}{2}}(-1)^j\binom{n}{2j+1}
S_{n,\;0}((n-1)^{2j}\lfloor\frac{x}{(n-1)^{n-1}}\rfloor))=$$
\begin{equation}\label{56}
\begin{cases} 0,\;\; if \; \lfloor\frac{x}{(n-1)^{n-1}}\rfloor \;\; is\;\; even,\\
1, \;\; if \;\; \lfloor\frac{x}{(n-1)^{n-1}}\rfloor \;\; is\;\; odd.\end{cases}
\end{equation}
Denote the left hand side of (\ref{56}) by $A_n(x).$ Then, similar to (\ref{45}),
we have
$$A_n(x)=$$
\begin{equation}\label{57}
 \begin{cases} 0,\;\; if \; x\equiv0,...,(n-1)^{n-1}-1,\pmod{2(n-1)^{n-1}},\\
1, \;\; if \;\;x\equiv (n-1)^{n-1},...,2(n-1)^{n-1}-1,\pmod{2(n-1)^{n-1}}.\end{cases}
\end{equation}
Furthermore, we consider the difference
 \begin{equation}\label{58}
\Delta_n(x)=S_{n,\;0}(x)-S_{n,\;0}((n-1)^{n-1}\lfloor\frac{x}{(n-1)^{n-1}}\rfloor).
   \end{equation}
   \begin{lemma}\label{l6}
   $(-1)^{s_{n-1}(x)}\Delta_n(x)$ is periodic with period $n(n-1)^{n-1}.$
   \end{lemma}
   \slshape Proof.\;\;\upshape Indeed, let
   $$x=n(n-1)^{n-1}t+j,\; j=0,1,...,n(n-1)^{n-1}-1.$$
   Let $j$ such that
   $$\lfloor\frac{j}{(n-1)^{n-1}}\rfloor=m, \; 0\leq m\leq n-1.$$
   Then
   $$j=(n-1)^{n-1}m+k,\; 0\leq k\leq (n-1)^{n-1}-1.  $$
   We have
$$\Delta_n(x)=S_{n,\;0}(n(n-1)^{n-1}t+j)-S_{n,\;0}(n(n-1)^{n-1}t+(n-1)^{n-1}m)=$$
$$S_{n,\;0}(n(n-1)^{n-1}t+(n-1)^{n-1}m+k)-S_{n,\;0}(n(n-1)^{n-1}t+(n-1)^{n-1}m)=$$

 \begin{equation}\label{59}
\sum_{i: (n-1)^{n-1}m+1\leq5i\leq (n-1)^{n-1}m+k-1}
(-1)^{s_4(n(n-1)^{n-1}t+5i)}.
\end{equation}

Note that
$$5i=(n-1)^{n-1}m+l,\;1\leq l\leq k-1\leq (n-1)^{n-1}-2.$$
Therefore, the summands in (\ref{59}) multiplied by $(-1)^{s_{n-1}(x)}$ have the
 form
 \newpage
$$(-1)^{s_{n-1}(n(n-1)^{n-1}t+(n-1)^{n-1}m+k)+
s_{n-1}(n(n-1)^{n-1}t+(n-1)^{n-1}m+l)}$$
and, since $l<k\leq (n-1)^{n-1}-1,$ this equal
$$(-1)^{s_{n-1}(n(n-1)^{n-1}t+(n-1)^{n-1}m)+s_{n-1}(k)
+s_{n-1}(n(n-1)^{n-1}t+(n-1)^{n-1}m)+s_{n-1}(l)}=$$
$$(-1)^{s_{n-1}(k)+s_{n-1}(l)}. $$
Therefore, the summands of (\ref{59}) not depend on $t$ and thus the sum (\ref{59}),
i.e., $\Delta_n(x)$ not depends on $t.$ \;$\blacksquare$
\begin{lemma}\label{l7}
The sequence
\begin{equation}\label{60}
\{(-1)^{s_{n-1}(x)+s_{n-1}(\lfloor\frac {x}{(n-1)^{n-1}}\rfloor-1)}A_n(x)\}
\end{equation}
is periodic with period $2(n-1)^{n-1}.$
   \end{lemma}
   \slshape Proof.\;\;\upshape In cases $x\equiv i\pmod{2(n-1)^{n-1}},\;
    i=0,1,...,(n-1)^{n-1}-1$ the terms of the
sequence are zeros. If $x\equiv i\pmod{2(n-1)^{n-1}},\;
i=(n-1)^{n-1},...,2(n-1)^{n-1}-1,$ put $x=2(n-1)^{n-1}t+i.$ Then $A_n(x)=1$
and we have
$$(-1)^{s_{n-1}(x)+s_{n-1}(\lfloor\frac {x}{(n-1)^{n-1}}\rfloor-1)}=(-1)^{s_{n-1}(2(n-1)^{n-1}t+i)+s_{n-1}(2t)}=$$

 $$(-1)^{s_{n-1}(2(n-1)^{n-1}t+i)+s_{n-1}(2(n-1)^{n-1}t)}=(-1)^{s_{n-1}(i)}$$
 and the lemma follows. \;$\blacksquare$\newline
 Now we obtain the following result.
\begin{theorem}\label{t8} If the conjectural relation (\ref{15}) is true, then we
have
\begin{equation}\label{61}
 S_{n,\;0}(x)=\sum_{j=0}^{\frac{n-3}{2}}(-1)^{\frac{n-3}{2}-j}\binom{n}{2j+1}
S_{n,\;0}((n-1)^{2j}\lfloor\frac{x}{(n-1)^{n-1}}\rfloor)+\nu_n(x),
\end{equation}
where $\nu_n(x)$ multiplied by $(-1)^{s_{n-1}(x)}$ is periodic with period
 $2n(n-1)^{n-1}.$
  \end{theorem}
\slshape Proof.\;\;\upshape Indeed, by (\ref{56})-(\ref{58}), we obtain (\ref{61})
with
$$\nu_n(x)=\Delta_n(x)
+(-1)^{\frac{n-1}{2}+s_{n-1}(\lfloor\frac{x}{(n-1)^{n-1}}\rfloor-1)}A_n(x).$$
Then, by Lemmas \ref{l6}-\ref{l7}, $(-1)^{s_{n-1}(x)}\nu_n(x)$ is periodic with
period equal the least common multiple of numbers $2(n-1)^{n-1}$ and
$n(n-1)^{n-1}.\; \blacksquare$\newline
As a corollary, in the case $n=3$ we again obtain Theorem \ref{t3} for
 $\nu(x)=\nu_3(x)$ but without detailed representation of $\Delta_3(x)$ and $\nu(x).$
\begin{remark}
It follows from the proof that, if for some
$$j=j_i,\; i=1,...,k,\; 1\leq j_1<j_2<...<j_k\leq\frac{n-3}{2},$$
to replace in (\ref{61}) $S_{n,\;0}((n-1)^{2j}\lfloor\frac{x}{(n-1)^{n-1}}
\rfloor) $ by $S_{n,\;0}(\lfloor\frac{x}{(n-1)^{n-1-2j}}\rfloor)$ and to
\newpage
denote
 the new sum by $\Sigma(j_1,...,j_k),$ then also the following form of
 Theorem \ref{t8} is valid
\begin{theorem}\label{t10} If the conjectural relation (\ref{15}) is true, then we
have
\begin{equation}\label{62}
 S_{n,\;0}(x)=\Sigma(j_1,...,j_k)+\nu_n^{(j_1,...,j_k)}(x),
\end{equation}
where $\nu_n^{(j_1,...,j_k)}(x)$ multiplied by $(-1)^{s_{n-1}(x)}$ is periodic
with period $2n(n-1)^{n-1}.$
  \end{theorem}
Thus we have $2^{\frac{n-3}{2}}$ different formulas of type (\ref{62}). In particular,
in case $n=3$ we have only formula, in case $n=5$ we have two different formulas, etc.
\end{remark}

 \section{Application of Theorem 8 in case $n=5$}
  Since the conjectural identity (\ref{15}) was proved in case $n=5,$ then, by
   Theorem \ref{t8}, we conclude that
\begin{equation}\label{63}
(-1)^{s_{4}(x)}\nu_5(x)=(-1)^{s_{4}(x)}(S_{5,\;0}(x)
-10S_{5,\;0}(16\lfloor
\frac{x}{256}\rfloor)+5S_{5,\;0}(\lfloor\frac{x}{256}\rfloor))
\end{equation}
is periodic with period 2560. If to write the period, then (\ref{63}) gives
 a recursion for $S_{5,\;0}(x).$ The computer calculations show that the period 
 with positions $\{0,...,2559\}$ contains all numbers from interval $[-35,35].$ Here
 we give several sequences of positions in $[0,2559]$ with these numbers 
 $g\in [-35,35]$.
 $$g=-35:\{251,252,254\},$$
 $$g=-34:\{246,249,1531,1532,1534\},$$
 $$g=-33:\{241,243,244,1526,1529\},$$
 $$g=-32:\{237,239,1521,1523,1524\},$$
 $$g=-31:\{231,232,234,1517,1519\},$$
 $$g=-30:\{197,199,200,217,219,220,226,229,511,1511,1512,1514,$$
 $$2497,2499,2500,2557,2559\},$$
 $$...$$
 $$g=30:\{196,198,216,218,227,228,230,1513,$$
 $$1515,2496,2498,2556,2558\},$$
 $$g=31:\{233,235,1516,1518,1520\},$$
 $$g=32:\{236,238,240,1522,1525\},$$
 $$ g=33:\{242,245,1527,1528,1530\},$$
 $$ g=34:\{247,248,250,1533,1535\},$$
 $$g=35:\{253,255\}.$$
Besides, by Theorem \ref{t10}, also
\newpage
\begin{equation}\label{64}
(-1)^{s_{4}(x)}\nu_5^{(1)}(x)=(-1)^{s_{4}(x)}(S_{5,\;0}(x)
-10S_{5,\;0}(\lfloor
\frac{x}{16}\rfloor)+5S_{5,\;0}(\lfloor\frac{x}{256}\rfloor))
\end{equation}
is periodic with period 2560. Again, if to write the period, then (\ref{64}) gives
 another recursion for $S_{5,\;0}(x).$ The computer calculations show that the 
 period with positions $\{0,...,2559\}$ contains all numbers from interval
  $[-9,9].$ Several sequences of positions in $[0,2559]$ with these numbers
   $h\in [-9,9]$ are the following:
 $$ h=-9: \{2411,2412,2414,2491,2492,2494\},$$
 $$ h=-8: \{1131,1132,1134,1211,1212,1214,2406,2409,2486,2489\},$$
 $$...$$
 $$ h=8: \{1133,1135,1213,1215,2407,2408,2410,2487,2488,2490\},$$
 $$ h=9: \{2413,2415,2493,2495\}.$$
 Finally, note that the sequence of the numbers of different values of\newline $\nu_3(x),\;
 \nu_5^{(1)}(x), \;\nu_5(x), etc.$ begins with $\{5,19,71,...\}\;.$
 \section{Recursions for $S_{3,\;1}(x)$ and $S_{3,\;2}(x)$}
Using (\ref{22})-(\ref{24}), it is easy to show that the form $3y(x)-y(4x)$ is
  invariant with respect to $S_{3,\;i}(x),\;i=0,1,2.$ This means that together with
 \begin{equation}\label{65}
 3S_{3,\;0}(x)-S_{3,\;0}(4x)=S_{3,\;0}(x)+S_{3,\;1}(x)+S_{3,\;2}(x),
  \end{equation}
  we have also
 \begin{equation}\label{66}
 3S_{3,\;1}(x)-S_{3,\;1}(4x)=S_{3,\;0}(x)+S_{3,\;1}(x)+S_{3,\;2}(x),
  \end{equation}
  \begin{equation}\label{67}
 3S_{3,\;2}(x)-S_{3,\;2}(4x)=S_{3,\;0}(x)+S_{3,\;1}(x)+S_{3,\;2}(x).
  \end{equation}
  Using (\ref{66})-(\ref{67}), as in Section 6, we can prove that the expressions
 \begin{equation}\label{68}
  (-1)^{s_2(x)}(S_{3,\;1}(x)-3S_{3,\;1}(\lfloor\frac{x}{4}\rfloor)),
   \end{equation}
   and
   \begin{equation}\label{69}
  (-1)^{s_2(x)}(S_{3,\;2}(x)-3S_{3,\;2}(\lfloor\frac{x}{4}\rfloor)),
 \end{equation}
  are eventually priodic with the same period as $(-1)^{s_2(x)}\nu(x)$ (\ref{12}),
  i.e., the period (\ref{55}), such that for $S_{3,\;2}(x)$ the period starts at $x=8,$ while for $S_{3,\;1}(x)$ the period starts at $x=16.$
  This means that, for $S_{3,\;i}(x),\;i=1,2,$ the same recursions hold as the recursion
  for $S_{3,\;0}(x)$ (\ref{11}) with the same function $\nu(x)$ (\ref{12}):
 \begin{equation}\label{70}
S_{3,1}(x)=3S_{3,1}\left(\left\lfloor \frac x
4\right\rfloor\right)+\nu(x),\;x\geq16,
\end{equation}
\newpage
with the initials
\begin{equation}\label{71}
S_{3,1}(x)=\begin{cases} 0,\;\; if \; x=0,1,\\
-1, \;\; if \; x=2,3,4, \\-2, \;\; if
 \; x=5,6,7,11,12,13,\\-3,\;\; if
 \; x=8,9,10,14,15. \end{cases}
 \end{equation}
 \begin{equation}\label{72}
S_{3,2}(x)=3S_{3,2}\left(\left\lfloor \frac x
4\right\rfloor\right)+\nu(x),\;x\geq8,
\end{equation}
with the initials
\begin{equation}\label{73}
S_{3,2}(x)=\begin{cases} 0,\;\; if \; x=0,1,2,6,7,\\
-1, \;\; if \; x=3,4,5. \end{cases}
 \end{equation}
For example, by (\ref{70}), (\ref{71}) and (\ref{12}), we have
$$S_{3,1}(20)=3S_{3,1}(5)+\nu(20)=3\cdot(-2)+(-1)^{s_2(20)}=-5; $$
analogously, by (\ref{72}), (\ref{73}) and (\ref{12}), we find
$$S_{3,2}(20)=3S_{3,2}(5)+\nu(20)=3\cdot(-1)+(-1)^{s_2(20)}=-2. $$
\section{A generalization}
A generalization of the conjectural equality (\ref{15}) is the following
$$\sum_{j=0}^{\frac{n-1}{2}}(-1)^j\binom{n}{2j+1}S_{n,\;i}((n-1)^{2j}x)=$$
\begin{equation}\label{74}
\sum_{j=0}^{n-1}S_{n,\;j}(x),\; i=0,...,n-1,\;x\geq1,\; n\equiv1\pmod2.
\end{equation}
If this conjecture is valid, then, as in the previous sections, we can obtain
the same recursions for every digit function $S_{n,\;i}(x),\; i=1,...,n-1,$
as for $S_{n,\;0}(x)$ (cf. Theorems \ref{t8}, \ref{10}). The question on initials
in cases $i\geq1$ we here remain open.

\end{document}